# ENGINEERING SOLUTIONS FOR NON-STATIONARY GAS PIPELINE RECONSTRUCTION AND EMERGENCY MANAGEMENT

**I.G. Aliyev**

*Operation and Reconstruction of Buildings and Facilities Department, Azerbaijan University Architecture and Construction, Baku, Azerbaijan, i-q-aliyev@mail.ru*

**Abstract-** The reconstruction, management, and optimization of gas pipelines is of significant importance for solving modern engineering problems. This paper presents innovative methodologies aimed at the effective reconstruction of gas pipelines under unstable conditions. The research encompasses the application of machine learning and optimization algorithms, targeting the enhancement of system reliability and the optimization of interventions during emergencies. The findings of the study present engineering solutions aimed at addressing the challenges in real-world applications by comparing the performance of various algorithms. Consequently, this work contributes to the advancement of cutting-edge approaches in the field of engineering and opens new perspectives for future research. A highly reliable and efficient technological Figure has been proposed for managing emergency processes in gas transportation based on the principles of the reconstruction phase. For complex gas pipeline systems, new approaches have been investigated for the modernization of existing control process monitoring systems. These approaches are based on modern achievements in control theory and information technology, aiming to select emergency and technological modes. One of the pressing issues is to develop a method to minimize the transmission time of measured and controlled data on non-stationary flow parameters of gas networks to dispatcher control centers. Therefore, the reporting Figures obtained for creating a reliable information base for dispatcher centers using modern methods to efficiently manage the gas dynamic processes of non-stationary modes are of particular importance.

## 1. INTRODUCTION

The increasing consumption volumes and resource usage to meet the gas demand in the country can be ensured not only through the design and construction of new networks and structures but also through the reconstruction of gas pipelines and efficient management by utilizing potential reserves. Gas pipelines play a critical role in the energy landscape; however, their reconstruction and management are fraught with challenges. This paper thoroughly examines the engineering aspects of gas pipeline reconstruction, including optimization and management methods. The recent application of advanced technologies, particularly machine learning, facilitates the resolution of engineering problems. Gas pipeline systems are complex systems characterized by numerous elements, and managing them requires a large amount of information corresponding to the number of possible elements. Managing such complex gas pipelines is impossible without a systematic approach based on the combined consideration of complex concepts of control theory, such as system, information, goal-oriented, and feedback. New approaches to reconstruction methods aimed at ensuring a higher level of quality are required in managing these systems.

The theory of selecting justified parameters for the operation and reconstruction of main gas pipelines requires the calculation of non-stationary processes resulting from any emergency [19, 20]. This task should be considered one of the most pressing issues. The operating conditions of a reconstructed gas pipeline must ensure its ability to perform specified functions reliably for a fixed period. The main purpose of the report Figure is to prepare a database in accordance with the requirements for the reconstruction and optimization technology of gas pipelines. This database can be applied to modern smart systems for managing gas supply using various methods and processes. It is known that recently, machine learning and artificial intelligence systems have been present in various fields, including the management of gas supply systems.

A comprehensive analysis of the current state of machine learning and artificial intelligence in solving gas supply issues has been conducted by several researchers [1, 2, 11, 14, 21]. These studies also discuss various types of machine learning and artificial intelligence methods that can be used to process and interpret data in different areas of the oil and gas industry. The use of these modern technologies can facilitate the decision-making process in the management of gas supply system reconstructions. The novelty of the work lies in the development of reporting Figures for solving several important new problems related to the management and regulation of technological processes in gas transportation based on optimization methods and mathematical modeling. For this purpose, an algorithm for the operational management of gas transportation systems has been developed. This algorithm is based on the information obtained about the emergency state of pipelines by monitoring changes in the technological parameters of gas flow at the peripheral parts of the pipelines during emergency modes. This significantly enhances the reliability of gas supply management and the efficiency of transportation processes during the reconstruction phase of complex gas pipelines.

At the current stage, the issues we face differ fundamentally from those already addressed [4, 8, 9, 15, 16, 17, 18] in that they involve imposing specific requirements and time constraints during the development phase of the reconstruction concept for complex gas pipeline systems. This ensures the acquisition of structural, technical, and technological bases within the framework of future operation of complex gas pipeline systems, the clarification and elimination of technical and technological contradictions in the process of gas transportation system development and reconstruction, and the efficient management of the gas transportation system under modern conditions.

## 2. PROBLEM STATEMENT

The essence of the studied method is to justify a targeted technological variant for reconstruction by increasing the efficiency of gas transportation under non-stationary conditions. Explores innovative approaches to modernizing existing control process monitoring systems within complex gas pipeline networks. These approaches integrate advancements in control theory and information technology to facilitate the selection of optimal emergency and operational modes. The reconstruction of gas transportation systems is characterized by the reliability of multi-level information bases transmitted to the dispatcher control center for complexes of local management systems of individual technological processes and objects.

A highly reliable and efficient technological model is introduced to manage emergency processes in gas transportation. This model is designed according to the principles outlined in the reconstruction phase, aiming to increase system efficiency while minimizing response times during emergency situations. A highly reliable and efficient technological model has been proposed for managing emergency processes in gas transportation, based on the principles of the reconstruction phase. For complex gas pipeline systems, new approaches have been investigated for modernizing existing control process monitoring systems. These approaches are based on modern advancements in control theory and information technology, aimed at selecting emergency and technological modes. Figure 1 shows methodological process for technological and emergency mode management in the gas pipeline systems.

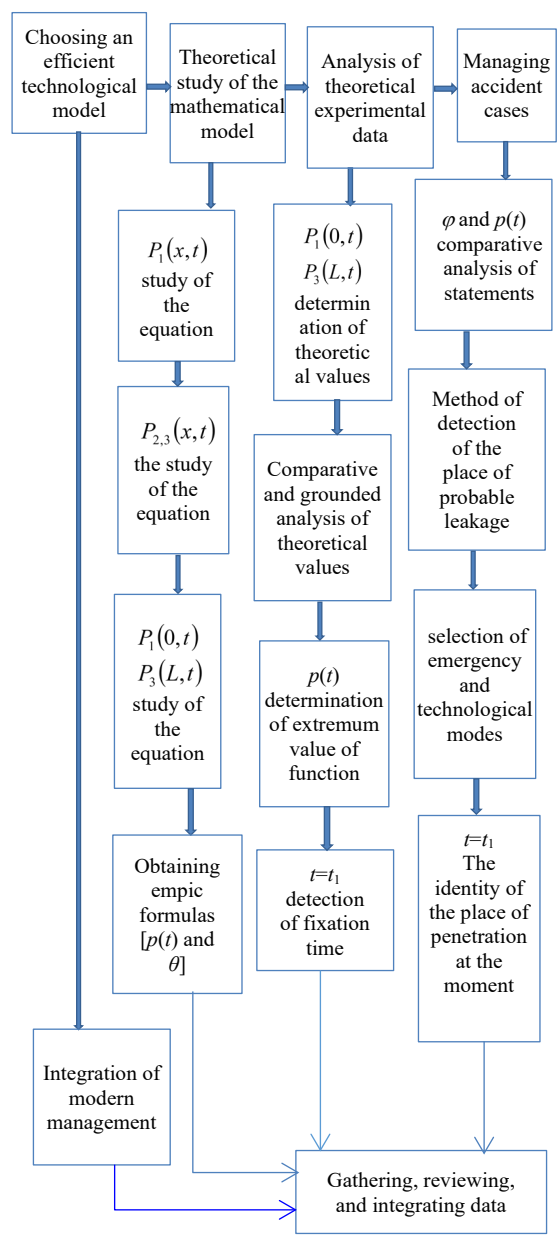

Figure 1. Process flow for gas pipeline management

Explores innovative approaches to modernizing existing control process monitoring systems within complex gas pipeline networks. These approaches integrate advancements in control theory and information technology to facilitate the selection of optimal emergency and operational modes. The reconstruction of gas transportation systems is characterized by the reliability of multi-level information bases transmitted to the dispatcher control center for complexes of local management systems of individual technological processes and objects. A major challenge identified in this study is the minimization of transmission time for measured and controlled data related to non-stationary flow parameters in gas networks. Efficient transmission of real-time data is essential for maintaining system integrity and safety during operation.

The reconstruction phase is based on the following key principles:
- The advantages of an efficient technological model for complex gas pipeline systems.
- The uninterrupted transportation of planned gas flows to consumers.
- The level of reliability indicators for complex gas pipeline systems.
- The safety of the operation and environmental aspects of complex gas pipeline systems.
- Leakage detection technology for pipeline safety.
- Integration of Modern Control Theories.

The "Main Pipeline Design" standards state that parallel gas pipelines are complex systems consisting of two or more pipeline lines. These lines have the same pressure at the beginning and end sections and are equipped with at least two connectors and connecting fittings to ensure the transfer of gas flow between the pipelines during technical operations and emergencies. In other words, the information obtained from changes in gas flow parameters at the end sections of the system due to the failure of one of the pipeline lines in parallel gas pipelines is not reliable. According to the research conducted in reference [5, 7], it was confirmed that depending on the operating conditions of gas pipelines, the pressure values at the beginning and end sections of damaged (ruptured) and undamaged pipeline lines were equal.

Moreover, it was even observed that, at any given time during operation, the change in pressure at the beginning and end of the gas pipeline was greater for the undamaged pipeline line than for the damaged one. In such conditions, identifying the damaged section in the dispatcher center becomes impossible. Indeed, the information obtained from the non-stationary nature of gas dynamic processes makes it unreliable for the effective management of gas pipeline systems. Therefore, the development of effective gas transportation Figures for the reconstruction of existing parallel gas pipelines has become a pressing issue. This is crucial for ensuring the long-term utilization and environmental safety of complex gas pipeline transportation, as well as for achieving reliable and trouble-free operation and efficient management of emergency situations.

To prepare the technological-based reporting Figure for this purpose, the following issues need to be addressed:
1) Implementation of effective strategies for managing parallel laid distributor pipeline systems during normal technological operating regimes and in case of emergencies;
2) Utilization of new methods with high quality for analyzing and optimizing the management of emergency and technological conditions and the selection of rational operating regimes;
3) Prevention and minimization of the consequences of emergency situations;
4) Ensuring the accuracy of damage location determination, automatic detection in real-time mode, and minimizing the impact of interventions in detection modes.

One of the pressing issues is to develop a method to minimize the transmission time of measured and controlled data on non-stationary flow parameters of gas networks to dispatcher control centers. The essence of the studied method lies in justifying a targeted technological variant for reconstruction by enhancing the efficiency of gas transportation under non-stationary conditions.

It is known that one of the main requirements for the reconstruction of gas pipelines is the application of modern technologies to the existing system. The China National Petroleum Corporation (CNPC) has independently developed 27 unique technologies across five series that have been successfully applied in the construction and management of various gas pipelines. These technologies have been implemented in the Lanzhou-Chengdu Chongqing Gas Pipeline, as well as pipeline projects in Sudan, Libya, India, Russia, and Central Asia. Some of these technologies include:
Oil and Gas Pipeline Flow Assurance Technology;
Oil and Gas Pipeline Simulation and Optimization Technology;
Pipeline Integrity Management System;
Early-warning and leakage detection technology for pipeline safety;
Pipeline Emergency Repair Technology, and so on.

However, leaks in pipeline networks remain one of the primary causes of countless losses to pipeline stakeholders and the surrounding environment. Pipeline failures can lead to severe environmental disasters, human casualties, and financial losses. To prevent such hazards and maintain a safe and reliable pipeline infrastructure, it is essential to implement early warning and leakage detection technology for pipeline safety. The pipeline leakage detection technologies are discussed in [6], summarizing the latest advancements. Various leak detection and localization methods in pipeline systems are reviewed, with their strengths and weaknesses highlighted. Comparative performance analysis is carried out to provide guidance on which leak detection method is most suitable for specific operational parameters. Furthermore, research gaps and open issues for the development of reliable pipeline leakage detection systems are discussed.

As an example, in [10, 18], we can mention a method based on a device that uses ultrasonic waves for radiation to manage pipelines. In the study being examined, acoustic sensors installed outside the pipeline detect the internal noise level and establish a baseline with specific characteristics. The self-similarity of this signal is continuously analyzed by the acoustic sensors. When a leak occurs, the resulting low-frequency acoustic signal is detected and investigated. If the properties of this signal differ from the original, an emergency signal is triggered. The received signal is stronger near the leak location, allowing for the localization of the leak. When a leak is detected, waves are reflected, and the location of the leak is determined by the emission speed of the wave and the return time of the scattered signals.

Ripple conversion is used to process these signals and locate the leak. In general, this method is based on the detection of noise generated in the pipeline when a leak occurs [18].

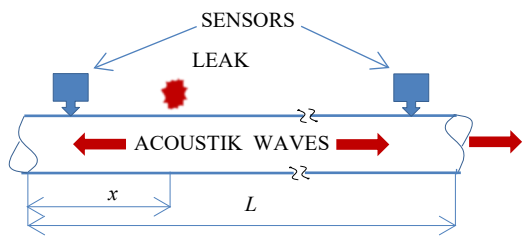

Figure 2. Technological figure of acoustic leak detection

Sensors are activated on the right and left of the point where the leak occurs, at any point $x$ on a pipeline with length $L$, as shown in Figure 2. In this study [10], technologies for extracting the characteristics of acoustic signals have been summarized. Experiments highlighted the characteristics in the time, frequency, and frequency-time domains by measuring acoustic leaks and interference signals.

Results indicate that external interference can be effectively eliminated through time domain, frequency domain, and frequency-time domain characteristics. It can be concluded that the acoustic leak detection method can be applied to natural gas pipelines and that its features can help reduce false alarms and missed alarms. Therefore, it can be concluded that the acoustic leak detection method can be applied to natural gas pipelines, helping to reduce false and missed alarm signals. Recently, many authors have extensively studied a system using an Internet of Things (IoT) analytical platform service to simulate real-time monitoring of pipelines and locate leaks within the pipeline [3]. In this study, pressure pulses based on the vibration principle of pipes are used for pipeline monitoring. The principle of time delay between pulse receptions at sensor positions is also applied in this research. The primary principle of the study was to create a wireless communication device that combines, programs, and interacts with Thing Speak IoT analytics platform through an Arduino & a Wi-Fi module.

A total of five channels were created to collect data from five sensors used in experimental testing on the platform (Figure 3). These signal data are collected every 15 seconds, and all channels are updated every 2 minutes. Thing Speak provides immediate visualization of data transmitted by the wireless communication device. Online analysis and processing of data were performed as it was received [3]. Using a data logging device and air as the transport fluid, the measured propagation speed of vibration was 355 m/s, ranging from 4.243 m to 4.246 m. This yielded only a 3 mm deviation when comparing the estimated leak location with the actual measured location of 4.23 m. The experimental results demonstrated that the performance of the wireless communication device compares favorably with that of the data logger and can determine the location of damage in actual pipelines when used for real-time monitoring.

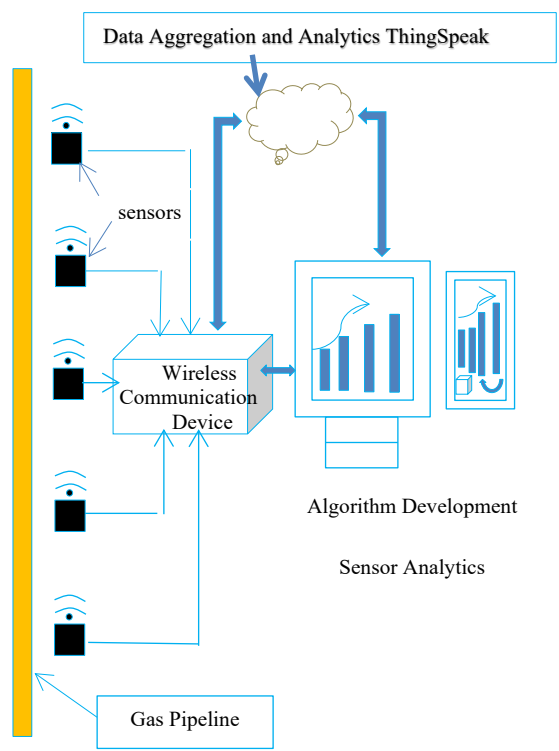

Figure 3. IoT system used for pipeline monitoring

With this communication device and analytical platform, it is possible to monitor pipelines in real time from anywhere in the world on any device with internet access. Integrity management programs are also available today that can locate the damage. A thorough, organized, and unified integrity management program me provides the means to improve the safety of pipeline systems. Such a program me provides the information for a pipeline system operator to allocate resources for prevention, detection, and mitigation actions that will result in improved safety and a reduction in the number of incidents as Figure 4 [22].

It covers performance and prescriptive based integrity management The performance-based integrity management program me employs more data and broad risk analysis, which permit pipeline system operators to comply with the requirements of this program me in the areas of inspection intervals, tools used, and mitigation techniques used. Inspection, prevention, detection, and mitigation are all part of the prescriptive process necessary to produce an integrity management program me.

Prior to proceeding with the performance-based integrity program me a pipeline system operator shall first ensure sufficient inspections are performed that provide the information on the pipeline condition required by the prescriptive based program me. The level of assurance of a performance-based program me shall meet or exceed that of a prescriptive program me [22]. It is known that there are three methods for locating leaks in gas pipelines: qualitative, quantitative, and analytical.

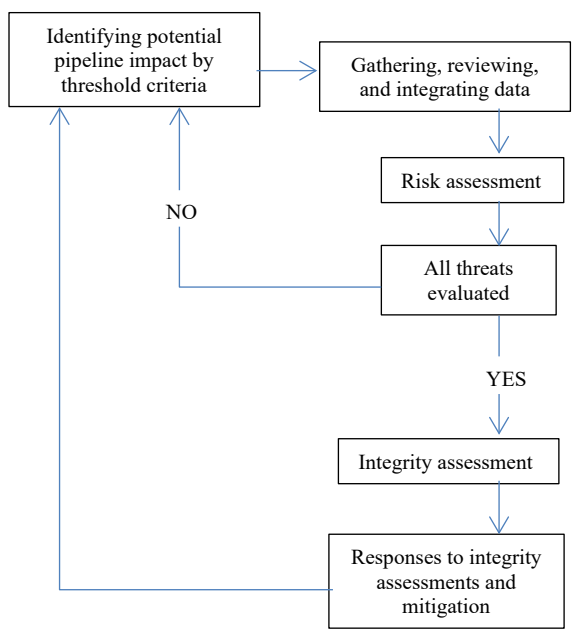

Figure 4. Integrity management plan process flow diagram

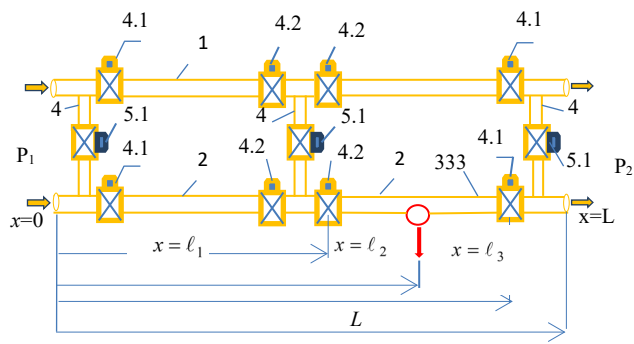

Figure 5. Principle diagram of a complex gas pipeline during the reconstruction phase

In Figure 5: 1- Undamaged pipeline ( $0 \leq x \leq L$ ); 2- Section from the starting point to the leakage point of the damaged pipeline ( $0 \leq x \leq \ell_2$ ), 3- Section from the leakage point to the end point of the damaged pipeline ( $\ell_2 \leq x \leq L$ ); 4- Connectors; 4.1 and 4.2 - Automatic valves on the pipeline; 5.1 and 5.2 - Automatic valves on the connectors; $P_1$ is pressure at the starting section of the gas pipeline, Pa; $P_2$ is pressure at the end section of the gas pipeline, Pa; $L$ is length of the gas pipeline, m, and $\ell_1, \ell_3$ are distance from the starting point of the damaged pipeline to the nearest shut-off fittings to the left and right of the leakage point, respectively.

In the proposed gas transportation Figure, at the fixed moment $t=t_1$, technological measures are implemented to isolate the damaged section from the main part of the pipeline. Specifically, the location of the damaged section is detected in the shortest possible time ( $\ell_2$ ), and the positions of the valves installed to the right ( $\ell_1$ ) and left ( $\ell_3$ ) of the leakage point are identified and closed by the dispatcher center at that moment. Based on the principle of managing technological measures, it is recommended to install automatic valves at the starting section of the pipeline network to the right of the first connector (4.1) and to the left of the last connector (4.2). In the normal (stationary) operating regime of gas pipelines, the valves on the connectors are closed, while the valves on the pipelines are open. In emergency situations (when the pipeline's integrity is compromised), the valves to the right and left of the leakage point on the pipeline are closed $(\ell_1, \ell_3)$, isolating the damaged section from the main part of the gas pipeline.

As shown in Figure 5 to ensure an adequate gas supply to consumers, the valves (5.1 and 5.2) located to the right and left of the leakage point on the connecting pipeline are opened under the control of dispatcher center management panels. The main distinguishing feature of the proposed gas pipeline compared to the existing one is its operation in an independent hydraulic regime. Another key difference of the reconstructed gas pipeline is the independence of the pipeline dynamics from each other during emergencies and technological regimes. The proposed parallel gas pipeline is a complex system characterized by a compressor station serving each pipeline at the starting section and numerous elements such as connectors and valves.

The practical examples mentioned mainly pertain to the quantitative method. In our research, however, the focus was on the analytical method. Specifically, it involved developing a reporting Figure to establish a reliable information base for dispatch centers to effectively manage the gas dynamic processes in non-stationary modes. As mentioned above, accidents occurring in parallel gas pipelines operating in a unified hydraulic regime cause a decrease in pressure throughout all sections of both pipeline networks over time. The decrease in pressure at the end section of the pipeline deteriorates the process of gas supply to consumers.

If we assume that the pressure difference is fully utilized by all consumers, then as the consumption of the initial gas mass in the network increases, the pressure at its end point will decrease again, leading to various devices and appliances experiencing different levels of disruption. Devices and appliances operating at maximum pressure will immediately sense the drop in gas pressure, reducing their numbers in operation simultaneously. In such a regime, only devices and appliances close to the nominal pressure will be able to operate. Therefore, these issues should be addressed in the reconstruction of existing gas pipelines. Based on the principles of the reconstruction phase, selecting an efficient technological figure that perfectly manages the technological and emergency processes of gas transportation is essential.

As we know, the reliability indicator of main gas pipelines directly depends on providing consumers with the required uninterrupted gas supply [5, 13]. Among the complex gas pipeline systems that excel at fulfilling this function is the parallel gas pipeline system. As a result of accidents, we adopt the following technological Figure for the efficient management of gas pipelines, the detection of leakage points, and the provision of uninterrupted gas supply to consumers in complex gas pipeline systems [7].

## 3. ESTABLISHING AN INFORMATION DATABASE FOR EFFICIENT MANAGEMENT OF THE RECONSTRUCTED GAS PIPELINE SYSTEM

During the reconstruction phase of complex gas pipelines, the acquisition of analytical expressions for operational management in an engineering context, along with the computational sequence for practical examples, has been investigated. Initially, known mathematical expressions reflecting the gas-dynamic processes of the transported gas flow under unsteady conditions were theoretically examined to justify the practical application of the proposed efficient technological model for the reconstruction of complex gas pipeline systems.

To substantiate the targeted technological variant for reconstruction, the process of pressure variation with respect to time along the pipeline length was studied in three specific sections (Figure 5). Two gas pipeline systems, with different parameters and suitable for practical operation, were selected for analysis. For the goal-oriented solution of the defined problems, simulations were conducted for various failure scenarios, where gas leakage occurred at different locations (near the start, end, and central parts) of the pipeline.

The pressure variation process in the beginning and end sections of both pipelines during emergency conditions was analyzed in detail. The results of parameter calculations for each pipeline variant were compared. The analysis revealed that despite differences in the dimensions and parameters of the selected gas pipelines, the results of the obtained analytical function and calculation performance remained consistent.

This finding demonstrates that regardless of the location of gas leakage during an emergency, the defined analytical function reaches its maximum and minimum values simultaneously in both pipelines. This indicator confirms the accuracy of the proposed analytical expressions and supports the practical applicability of the obtained data.

The fixed time duration is deemed acceptable for pipeline management at the agreed operational level. The analysis of the calculated values of the analytical expressions at the specified fixed time reveals that these values constitute a reliable and effective information base for the dispatcher control center. Thus, the analysis of the parameters at a particular time, $t=t_1$, enables the identification of whether the pressure changes in the end sections of the pipeline are due to an emergency or a change in the technological regime. In this context, the dispatcher control center can make an informed decision regarding the activation of the modern pipeline management system.

To establish an information database for the efficient management of the reconstructed gas pipeline system, the reporting Figure was implemented in the following sequence. First, let's assume that the integrity of the second pipeline in the parallel gas pipeline has been compromised, resulting in a large amount of gas leaking from the damaged section ($x = \ell_2$) into the environment (Figure 4).

At this time, the undamaged pipeline continues to operate in its previous stationary regime. In other words, the non-stationary regime that arises in the damaged pipeline does not affect the undamaged pipeline. As shown in Figure 4 the non-stationary regime of gas flow in the damaged gas pipeline due to the accident is divided into three different technological sections:

The section from the starting point of the pipeline to the shut-off valve to left of the leakage point ($0 \leq x \leq \ell_1$). The damaged section of the pipeline, between the shut-off valves to the left and right of the leakage point ($\ell_1 \leq x \leq \ell_3$). The section from the shut-off valve to the right of the leakage point to the end point of the pipeline ($\ell_3 \leq x \leq L$).

The hydraulic models of all three sections are different. The hydraulic model of the pipeline system is constructed from the models of its components: multi-line pipelines with loops and connectors, shut-off valves, compressor stations, etc. The models of these components, in turn, consist of the computational expressions for the elements. The solution to this problem necessitates considering the entire system as a whole rather than just a separate section of the pipeline. Initially, there is a significant need for analytical expressions to solve the equations describing non-stationary events in parallel gas pipelines and for mathematical expressions to solve these equations that reflect real events in the pipeline. The non-stationary operating conditions of pipelines lead to significant pressure changes and disruptions in the normal operation of gas delivery to consumers. Studying these processes to efficiently manage the pipelines and considering the results of the obtained analytical solution methods at dispatcher control points will significantly enhance the reliability of gas supply and the efficiency of transportation processes.

For this purpose, it is essential to obtain the analytic expressions reflecting the non-stationary events in the two-line parallel gas pipelines and ensure the management of physical processes through the solution of these equations. The mentioned analytic expressions have been determined based on the resolution of the problem addressed in our published work [7], which specifically deals with the non-stationary flow of gas in each of the three sections. These expressions are as follows [7].

$$P_1(x,t) = P_1 - 2aG_0 x + \frac{8aLG_0}{\pi^2} \sum_{n=1}^{\infty} \cos\frac{\pi nx}{L} \frac{e^{-(2n-1)^2 \alpha_3 t}}{(2n-1)^2} -$$
$$- \frac{4aL}{\pi^2} G_0 \sum_{n=1}^{\infty} ((1-(-1)^n) \cos\frac{\pi nx}{L} \times$$
$$\times \frac{e^{-n^2 \alpha_3 t}}{n^2} - 2aG_{ut}\left(\frac{x^2}{2L} + \frac{\ell_2^2}{2L} + \frac{L}{3} - \ell_2\right) +$$
$$+ \frac{4aL}{\pi^2} G_{ut} \sum_{n=1}^{\infty} \cos\frac{\pi nx}{L} \cos\frac{\pi n\ell_2}{L} \frac{e^{-n^2 \alpha_3 t}}{n^2} \quad (1)$$

$$P_{2,3}(x,t) = P_1(x,t) - 2aG_{ut}(x - \ell_2) \quad (2)$$

where, $\alpha_3 = \pi^2 c^2 / 2aL^2$.

By considering $x=0$ in Equation (1) and $x=L$ in Equation (2), we establish the conformity of the pressure distribution at the beginning and end points of the gas pipeline over time.

$$P_1(0,t) = P_1 - \frac{c^2 t}{L} G_{ut} + \frac{8aLG_0}{\pi^2} \sum_{n=1}^{\infty} \frac{e^{-(2n-1)^2 \alpha_3 t}}{(2n-1)^2} -$$
$$-2aG_{ut}\left(\frac{\ell_2^2}{2L} + \frac{L}{3} - \ell_2\right) + \tag{3}$$
$$+\frac{4aL}{\pi^2} G_{ut} \sum_{n=1}^{\infty} \cos\frac{\pi n \ell_2}{L} \frac{e^{-\alpha_2 t}}{n^2}$$

$$P_3(L,t) = P_2 - \frac{c^2 t}{L} G_{ut} + \frac{8aLG_0}{\pi^2} \sum_{n=1}^{\infty} (-1)^n \frac{e^{-(2n-1)^2 \alpha_3 t}}{(2n-1)^2} -$$
$$-2aG_{ut}\left(\frac{L}{2} + \frac{\ell_2^2}{2L} + \frac{L}{3} - \ell_2\right) + \tag{4}$$
$$+\frac{4aL}{\pi^2} G_{ut} \sum_{n=1}^{\infty} (-1)^n \cos\frac{\pi n \ell_2}{L} \frac{e^{-\alpha_2 t}}{n^2} - 2aG_{ut}(\ell_2 - L)$$

Based on Equations (3) and (4) and the results of the analysis conducted in [6], we obtain Equation (5):

$$\theta = \frac{1}{2} + \left[\frac{1-p(t)}{1+p(t)}\right] \times \left[\frac{2}{3} + \frac{e^{-2\alpha_2 t} - 4e^{-\alpha_2 t}}{\pi^2}\right] \tag{5}$$

where, $\alpha_2 = \frac{\pi^2 c^2}{2aL^2}$; $p(t) = \frac{P_1 - P_1(0,t)}{P_2 - P_2(L,t)}$; $\theta = \frac{\ell_2}{L}$

$G_{ut}$: Mass flow rate at the leakage points for the gas pipeline accident regime, Pa×sec/m

$G_0$: Mass flow rate at the beginning of the gas pipeline in the stationary regime, Pa×sec/m

$c$: Speed of sound propagation for gas in an isothermal process, m/s

$P$: Pressure at the ends of the gas pipeline, Pa

$x$: Coordinate along the length of the gas pipeline, m

$t$: Time coordinate, sec

$2a$: Charney linearization, 1/sec

It should be noted that the function $p(t)$ characterizes the ratio of pressure values at the beginning and end points of the gas pipeline. Based on the analysis of the parameters involved in Equation (5), it becomes evident that there is a functional dependence between $\theta$ and $p(t)$.

Based on the conducted research, it was determined that if $P(t) = 1$, then $\theta$ takes a value of 0.5. This indicates that the gas leakage occurs in the exact middle section of the pipeline. When $p(t)$ takes values greater than 1, $\theta$ takes values less than 0.5, confirming that the gas leakage occurs between the initial and middle sections of the pipeline. Conversely, when $p(t)$ takes values less than 1, $\theta$ takes values greater than 0.5, confirming that the gas leakage is situated between the middle and final sections of the pipeline.

To confirm the above ideas, the expressions (3) and (4) were examined, and the following parameters of the gas pipeline were considered to determine the law of pressure change along the gas axis at different characteristic points of gas leakage location (at the beginning, $\ell_2 = 0.5 \times 10^4$, m, in the middle, $\ell_2 = 5 \times 10^4$, m; and at the end, $\ell_2 = 9.5 \times 10^4$, m) [6]. $P_1 = 55 \times 10^4$ Pa; $P_2 = 25 \times 10^4$ Pa; $G_0 = 30$ Pa×sec/m; $c = 383.3 \frac{m}{\sec}$; $L = 10 \times 10^4$ m; $2a = 0.1$/sec.

Currently, the calculation of gas flow regimes is implemented through the use of modern software. Firstly, the pressure values of the gas pipeline at the beginning and end of the pipeline are calculated for each of the three characteristic points of the gas leakage every 100 seconds. The results are recorded in Table 1.

Table 1. Pressure values at the point of rupture of the gas pipeline depending on time

| ($t$) sec. | The location of the gas pipeline failure. | |
|---|---|---|
| | $\ell_2 = 0.5 \times 10^4$, m | |
| | The pressure in the starting part $P_1(0,t) \times 10^4$, Pa | Pressure in the last part $P_3(L,t) \times 10^4$, Pa |
| 1 | 2 | 3 |
| 100 | 52.23 | 25 |
| 200 | 50.58 | 25 |
| 300 | 49.3 | 24.99 |
| 400 | 48.21 | 24.97 |
| 500 | 47.25 | 24.93 |
| 600 | 46.38 | 24.84 |
| 700 | 45.58 | 24.72 |
| 800 | 44.83 | 24.57 |
| 900 | 44.13 | 24.37 |
| | $\ell_2 = 5 \times 10^4$, m | |
| 1 | 2 | 3 |
| 100 | 55 | 24.99 |
| 200 | 54.9 | 24.89 |
| 300 | 54.66 | 24.66 |
| 400 | 54.34 | 24.33 |
| 500 | 53.96 | 23.96 |
| 600 | 53.56 | 23.49 |
| 700 | 53.14 | 23.14 |
| 800 | 52.71 | 22.71 |
| 900 | 44.13 | 24.37 |
| | $\ell_2 = 9.5 \times 10^4$, m | |
| 1 | 2 | 3 |
| 100 | 55 | 22.23 |
| 200 | 55 | 20.58 |
| 300 | 54.99 | 19.3 |
| 400 | 54.97 | 18.21 |
| 500 | 54.925 | 17.25 |
| 600 | 54.84 | 16.38 |
| 700 | 54.72 | 15.58 |
| 800 | 54.565 | 14.84 |
| 900 | 54.37 | 14.14 |

Using the data provided in Table 1, we plot the Figure of $p(t)$ as a function of time for each of the three characteristics. The analysis of Figure 6 reveals that, regardless of the location of the gas pipeline rupture based on the parameters mentioned above, the $p(t)$ function attains maximum or minimum values within a duration of $t=300$ seconds. Therefore, the fixed time observed at dispatcher stations for the examined gas pipeline is $t_1 = 300$ seconds.

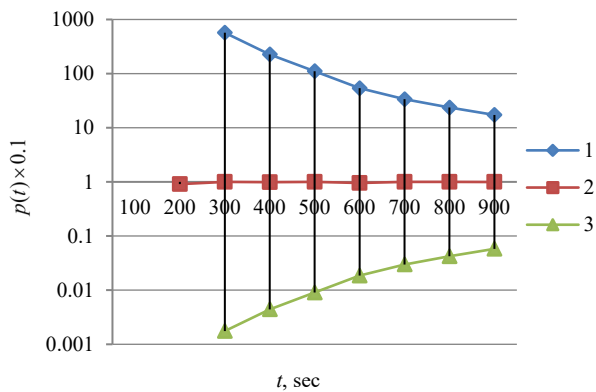

Figure 6. Variation pattern of *p(t)* function over time, dependent on the location of the gas pipeline rupture, (1- $\ell_2 = 0.5 \times 10^4$, m;
2- $\ell_2 = 5 \times 10^4$, m; 3- $\ell_2 = 9.5 \times 10^4$, m)

Minimizing the transmission time of information to dispatcher control centers is crucial. Therefore, the adopted value of $t_1$= 300 seconds ensures the safe operation of gas pipelines. On the other hand, the *p(t)* function takes different values depending on the dependence on $\theta$, which determines the location of the leakage. Specifically, when $\theta > 0.5$, the *p(t)* function takes its minimum value (3rd line in Figure 6), while for $\theta < 0.5$, the *p(t)* function takes its maximum value (1st line in Figure 6). However, when $\theta$ is approximately equal to 0.5, *p(t)* approaches 1, indicating a scenario where the leakage point is in the middle section of the pipeline ($\ell_2 = 5 \times 10^4$, m).

The analysis of the figure leads to the following conclusion: when *p(t)* belongs to the interval [0;1], it characterizes the location of gas pipeline rupture occurrence between $0 < \ell_2 < \frac{L}{2}$. When *p(t)* is greater than 1, it confirms the occurrence of gas pipeline rupture between $\frac{L}{2} < \ell_2 < L$.

When *p(t)* = 1, it confirms that the location of gas pipeline rupture is at the exact midpoint ($\ell_2 = \frac{L}{2}$).
Clearly, these simple and understandable pieces of information are very important for the gas emergency control center. Based on the analysis conducted, it was determined that the fixation period $t=t_1$ remains constant regardless of the location of accidents on the pipeline. To determine whether this time remains constant or varies depending on the parameters of the pipeline, we consider different parallel gas pipelines with the following parameters particular importance.
$P_1$=14×10$^4$ Pa; $P_2$ =11×10$^4$ Pa; $G_0$= 10 Pa×sec/m; $c = 383.3 \frac{m}{sec}$; $L$ =3×10$^4$ m

It is evident from the given information that the length of the newly adopted parallel gas pipeline is at least 3 times shorter than the length of the previously analyzed gas pipeline, and also that the values of other parameters (pressure and consumption) are different. We consider the length of the examined complex parallel pipeline at three characteristic points where the gas leakage occurs, similar to the previous pipeline: near the beginning of the pipeline $\ell_2 = 0.5 \times 10^4$, m at the middle $\ell_2 = 1.5 \times 10^4$, m and near the end $\ell_2 = 2.5 \times 10^4$, m.

We use Equations (3) and (4) to calculate the pressure values of the gas pipeline's beginning and end sections every 60 seconds for each of the three characteristic points where the gas leakage occurs. The results are recorded in Table 2.

Table 2. Pressure values of the gas pipeline depending on the location of the leakage and time

| *t* (sec) | The location of the gas pipeline failure | |
|---|---|---|
| | $\ell_2 = 0.5 \times 10^4$, m | |
| | The pressure in the starting part $P_1(0,t) \times 10^4$, Pa | Pressure in the last part $P_3(L,t) \times 10^4$, Pa |
| 1 | 2 | 3 |
| 60 | 13.37 | 10.97 |
| 120 | 12.95 | 10.79 |
| 180 | 12.61 | 10.55 |
| 240 | 12.29 | 10.27 |
| 300 | 11.99 | 9.98 |
| 360 | 11.70 | 9.69 |
| 420 | 11.40 | 9.40 |
| 480 | 11.11 | 9.11 |
| 540 | 10.81 | 8.81 |
| 600 | 10.52 | 8.52 |
| | $\ell_2 = 1.5 \times 10^4$, m | |
| 1 | 2 | 3 |
| 60 | 13.83 | 10.83 |
| 120 | 13.54 | 10.54 |
| 180 | 13.24 | 10.24 |
| 240 | 12.95 | 9.95 |
| 300 | 12.66 | 9.66 |
| 360 | 12.36 | 9.36 |
| 420 | 12.07 | 9.07 |
| 480 | 11.77 | 8.77 |
| 540 | 11.48 | 8.48 |
| 600 | 11.19 | 8.19 |
| | $\ell_2 = 2.5 \times 10^4$, m | |
| 1 | 2 | 3 |
| 60 | 13.97 | 10.37 |
| 120 | 13.79 | 9.95 |
| 180 | 13.55 | 9.61 |
| 240 | 13.27 | 9.29 |
| 300 | 12.98 | 8.99 |
| 360 | 12.69 | 8.70 |
| 420 | 12.40 | 8.40 |
| 480 | 12.11 | 8.11 |
| 540 | 11.81 | 7.81 |
| 600 | 11.52 | 7.52 |

Using the data from Table 2, we plot the variation of *p(t)* as a function of time for each of the three characteristic points (Figure 7). From the analysis of Figure 7, it can be noted that the variation of pressure over time in the outer sections depending on the location of the gas pipeline's leakage is similar to that observed in Figure 6.

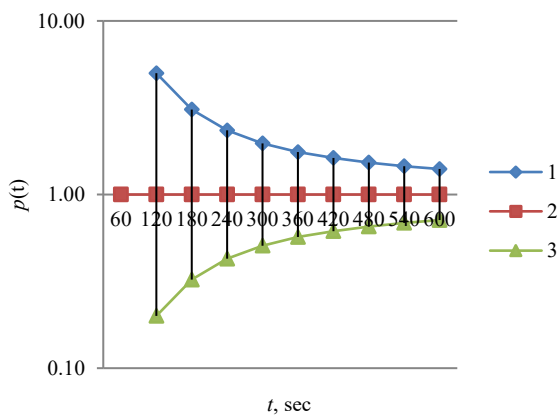

Figure 7. The variation pattern of $p(t)$ function depending on the time and the location of the gas pipeline's leakage
(1- $\ell_2 = 0.5\times10^4$, m ; 2- $\ell_2 = 15\times10^4$, m ; 3- $\ell_2 = 2.5\times10^4$, m )

From the analysis of Figure 7, it can be noted that the variation of pressure over time in the outer sections depending on the location of the gas pipeline's leakage is similar to that observed in Figure 6. Specifically, when $\theta > 0.5$, the $p(t)$ function takes its minimum value (indicated by the 3rd line in Figure 7), while for $\theta < 0.5$, the $p(t)$ function takes its maximum value (indicated by the 1st line in Figure 7).

However, regardless of the location of the gas pipeline's leakage, the $p(t)$ function takes its maximum or minimum value over a duration of $t=120$ seconds. This means that the fixed duration $t_1=120$ seconds is applicable for this gas pipeline. It can be concluded that the fixed duration varies depending on the length of the gas pipeline and the parameters of the gas flow. Thus, as the length of the pipeline segment decreases, this duration decreases accordingly. This is because the time required for the parameters of the non-stationary gas flow during the accident to reach the starting point of the gas pipeline is related to the change in the parameters. Considering the short length of the analyzed gas pipeline segment and the reliability of the information transmitted to the dispatcher station, it is necessary to accept a minimum duration of fixation that is sufficient to record the change in pressure of the non-stationary gas flow at the dispatcher station.

In other words, even if the speed of the non-stationary gas flow in the pipeline matches its sound propagation speed, the delivery time of information from the last three points of a 100000-meter pipeline segment to the starting point cannot be faster than 261 seconds. It should be noted that considering this indicator during the reconstruction stage of gas pipelines for efficient management allows for the optimization of the issue of minimizing the transfer time of information and control system data to dispatcher control centers.

We conduct theoretical research to examine the consistency of the method for determining the fixing time $t=t_1$, which is one of the measures taken to reduce the workload of the dispatcher station. For this purpose, we perform engineering calculations based on theoretical research. Using Equation (5), we determine the location of the pipeline's instability point at each of the three characteristic points, depending on time. First, we calculate based on the data of the gas pipelines under consideration, with a length of $L$=100000 m, and record the calculation results in Table 3.

Table 3. Relative error of the actual and calculated values of the gas instability point location depending on time for a pipeline with a length of 100,000 meter

| $t$ (sec) | Actual values of failure of the pipeline with a length of 100000 m | |
|---|---|---|
| | $\ell_2 = 0.5\times10^4$, m | |
| | Calculated by Equation (5) price | Relative error of actual and calculated prices |
| 1 | 2 | 3 |
| 100 | - | - |
| 200 | - | - |
| 300 | 0.55 | 0.04 |
| 400 | 0.33 | 3.44 |
| 500 | 0.15 | 1.61 |
| 600 | 0.04 | 1.38 |
| | $\ell_2 = 5\times10^4$, m | |
| 1 | 2 | 3 |
| 100 | 5.00 | 0.00 |
| 200 | 5.20 | 0.02 |
| 300 | 5.00 | 0.00 |
| 400 | 5.04 | 0.00 |
| 500 | 5.00 | 0.00 |
| 600 | 5.12 | 0.01 |
| | $\ell_2 = 9.5\times10^4$, m | |
| 1 | 2 | 3 |
| 100 | 5.00 | 0.47 |
| 200 | 9.20 | 0.03 |
| 300 | 9.45 | 0.01 |
| 400 | 9.67 | 0.02 |
| 500 | 9.84 | 0.04 |
| 600 | 9.96 | 0.05 |

Sequentially, based on the given data for the next gas pipeline with a length of $L$= 30000 m, we calculate the locations of instability points at each of the three characteristic points using Equation (5) and record the results in Table 4. Based on the analysis of Tables 3 and 4, it can be noted that the prescribed Equation (5) is useful for engineering calculations. The results of both analyzed gas pipelines' calculations showed that only at the fixed time $t=t_1$ is the location of instability accurately determined, with its relative error value not exceeding 0.04, or in other words, less than 4%.

Therefore, the prescribed calculation Figure is a suitable method for verifying the compliance of gas pipelines with the requirements for safe and efficient utilization, determining the parameters of gas flow in non-stationary regimes, as well as selecting appropriate parameters for the reconstruction of gas pipelines.

Table 4. The values of the relative error of the actual and calculated positions of instability of the gas pipeline with a length of 30,000 m, depending on time

| $t$ (sec.) | Actual values of failure of the pipeline with a length of 30000 m | |
|---|---|---|
| | $\ell_2 = 0.5 \times 10^4$, m | |
| | Calculated by Equation (5) | Relative error of actual and calculated prices |
| 1 | 2 | 3 |
| 60 | 0.06 | 0.87 |
| 120 | 0.28 | 0.44 |
| 180 | 0.51 | 0.02 |
| 240 | 0.71 | 0.41 |
| 300 | 0.85 | 0.70 |
| 360 | 0.95 | 0.91 |
| | $\ell_2 = 1.5 \times 10^4$, m | |
| 1 | 2 | 3 |
| 60 | 1.50 | 0.0 |
| 120 | 1.50 | 0.0 |
| 180 | 1.50 | 0.0 |
| 240 | 1.50 | 0.0 |
| 300 | 1.50 | 0.0 |
| 360 | 1.50 | 0.0 |
| | $\ell_2 = 2.5 \times 10^4$, m | |
| 1 | 2 | 3 |
| 60 | 2.94 | 0.175 |
| 120 | 2.72 | 0.088 |
| 180 | 2.49 | 0.005 |
| 240 | 2.29 | 0.083 |
| 300 | 2.15 | 0.140 |
| 360 | 2.05 | 0.181 |

## 4. DETERMINATION OF ANALYTICAL EXPRESSION BASED ON WHICH ACCIDENTS AND TECHNOLOGICAL REGIMES ARE SELECTED FOR DISPATCHER STATIONS

The variation in pressure in a gas pipeline can result from both changes in the amount of gas supplied to consumers, the opening and closing of valves in the pipeline, and the operation of various elements within the pipeline. As a result, technological processes cause changes in pressure at the beginning and end sections of the gas pipeline. One of the most important tasks of operational control in gas pipeline technological regimes is to identify critical conditions in gas supply systems at dispatcher stations. For instance, pressure variations at the end sections of the pipeline can occur as a result of changes in technological regimes. Therefore, it is essential to prioritize understanding the pressure changes resulting from technological processes in the event of accidents.

The initial model of events and processes occurring in gas transmission systems is not always accurately represented with high precision. Therefore, to solve the problems of operational management of gas transmission systems, it is necessary to use analytical expressions to assess the real situation of gas transmission systems. Having a mathematical expression allows for selecting the parameters and control structures of the pipeline, determining optimality criteria and constraints, assessing accuracy, selecting appropriate technical control devices, and so on.

To solve the problems of managing trunk gas pipelines, it is necessary to understand the non-stationary dynamic characteristics of these systems. It is appropriate to utilize the formalization of technological and accident processes for the transportation of gas in pipeline systems. By employing the method of describing accident processes, it is possible to obtain unit characteristics. Based on theoretical research, the following inequality has been determined.

$$\frac{\varphi - 0.5}{\varphi + 0.5} < p(t) < \frac{\varphi + 0.5}{\varphi - 0.5} \qquad (6)$$

where, $\varphi = \frac{2}{3} + \frac{e^{-2\alpha_2 t} - 4e^{-2\alpha_2 t}}{\pi^2}$.

It has been determined that if the function $p(t)$ satisfies the conditions of the inequality mentioned above, the variation of pressure over time corresponds to the lawfulness characterizing accidents in the gas pipeline, otherwise, it is the result of technological processes. Therefore, the function $p(t)$ is a reliable parameter that enables dispatchers to make timely decisions for the efficient management of the gas pipeline. At the point where time $t=t_1$ is fixed, the value of $\theta$\theta and the obtained values of the function $p(t)$ are reliable information for the dispatcher station. This is because the accuracy of determining the leakage location is characterized by the moment $t=t_1$. To visualize the methodology, we use the flow diagram below.

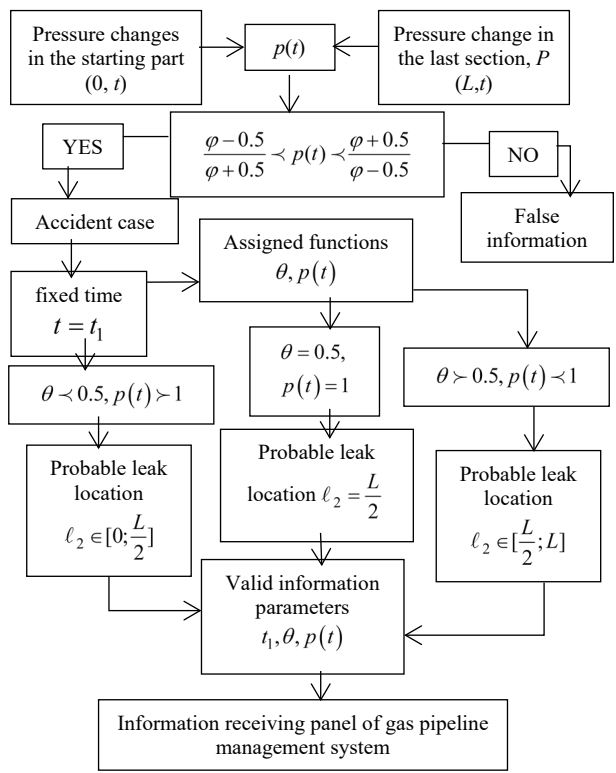

Figure 8. A visual diagram of the methodology

On the other hand, for the reliability and effectiveness of automatic detection of the accident location in online mode, these parameters are useful. For instance, if pressure changes are observed in the initial and final sections of the gas pipeline, the dispatcher station first determines whether this change occurred as a result of the accident based on the conditions of inequality (6). Sequentially, the values of the function $p(t)$ are automatically calculated, and its maximum or minimum value is determined.

The obtained values are assumed to be equal to $t=t_1$. Based on the value of $p(t)$ at $t=t_1$ calculated according to the inequality (5), the value of $\theta$ is determined. As shown in Figure 8. If this value is less than 0.5 and $p(t)>1$, it is determined that the location of the accident is between the starting and middle sections of the gas pipeline, and based on the value of $\theta$, the location of the gas leak $(\theta = \frac{\ell_2}{L})$ is automatically identified. Also, the valves $(\ell_1, \ell_3)$, whose locations are already known, are activated to isolate the damaged section of the pipeline from the main section of the gas pipeline.

Based on the results of the research, it can be stated that using expressions (5) and (6), it is possible to distinguish between accidental and technological regimes. On the other hand, from the moment the occurrence of an accident is confirmed, it is possible to determine the location of gas leaks and the sequence of valve activations. Thus, it can be concluded that the values obtained at the fixed time $t=t_1$ and from these expressions provide a reliable information base for monitoring and efficiently managing the technological regimes of gas distribution pipelines under extreme conditions.

To ensure the comprehensive management of system integrity, it is essential that emergency dispatch center personnel accurately understand whether pressure fluctuations in the final segments of the gas pipeline are due to incidents or technological processes and can discern the relationships between these processes. In this context, modern control methods, particularly those enhanced by artificial intelligence (AI) technologies, should be prioritized. For instance, AI technologies, particularly hybrid models, provide a more resilient and effective approach by offering accurate predictions even with limited data.

In this study, we compare the performance of algorithms by presenting tables that highlight the accuracy, efficiency, and other relevant metrics of each algorithm (Tables 3 and 4). Specifically, indicators that assess algorithms' capabilities in real-time incident detection and data augmentation for handling data insufficiencies have been presented. This approach lays a stronger foundation for addressing existing challenges in pipeline engineering.

Furthermore, key tasks include tracking the evolution of incidents, using data augmentation to address data gaps, enhancing model interpretability, and improving the management of pipeline integrity. Future research should focus on the optimization and integration of AI models to support informed decision-making based on a reliable information base. Additionally, AI enables the analysis of large volumes of data from various sources to predict incident risks, optimize management processes, and provide critical insights that support engineering decisions. Consequently, prioritizing the broader application of hybrid models in real-time monitoring and predictive technologies is essential for advancing this field.

## 5. CONCLUSIONS

A highly efficient technological Figure with a high reliability indicator has been proposed, based on the principles of the reconstruction phase, for the management of technological processes in gas transportation. The application of modern technologies and equipment to the proposed gas transportation Figure is of particular importance for optimizing the gasification capabilities of consumers.

Several mathematical models have been developed and an effective computational Figure has been devised to address a number of important new problems related to the management and regulation of technological processes in gas transportation. This will significantly enhance the reliability and efficiency of gas supply management during the reconstruction phase of complex gas pipelines. The calculation results showed that the error in comparing the obtained algorithm values does not exceed 4%. The developed computational Figure can be used to create an information database for the efficient management of gas transportation systems.

To ensure the unconditional resolution of efficient management issues, a function characterized by the ratio of pressure drops has been defined and tested. The analysis of the $p(t)$ function and parameters, which can perform multiple operations at the fixed time $t=t_1$, demonstrated that considering this indicator during the reconstruction phase can optimize the minimization of the data transmission time from information and control systems to dispatcher control points. The established relationship enables the differentiation of technological (stationary) and non-stationary operating modes of gas pipelines based on the indicators of changes in the parameters of the gas flow's final sections.

Based on the results of the analyses, individual problems related to the reconstruction of parallel main gas pipelines have been solved using the developed algorithms and methods, taking into account the non-stationary regime. Consequently, targeted methods for the reconstruction processes of an efficiently managed gas transport system have been identified as the subject of the research.

A technological Figure has been proposed that adopts the integration of new technologies and methodologies to ensure efficiency, safety, and reliability for the effective management of the reconstructed gas pipeline system. A method for the rapid detection of faults, which should become an integral part of the automated dispatch control system of the gas pipeline, is proposed. This method involves recording the gas flow pressure at the start and end of the pipeline at a fixed moment $t=t_1$.

The applied optimization algorithms not only enhance system reliability but also optimize interventions during emergencies. Future research could further advance by exploring the application of these methods in other engineering fields and integrating new technologies. This work offers new perspectives for improving the efficiency of gas pipeline systems and lays the groundwork for future investigations.

The proposed technological model provides a comprehensive solution for managing complex gas transportation systems. Future work should focus on further improving the model by integrating additional predictive technologies and enhancing data transmission methods. Furthermore, continuous advancements in control systems and emergency response strategies are necessary to address the increasing complexity of global gas transportation networks.

**REFERENCES**

[1] A.B. Sircar, K. Yadav, K. Rayavarapu, N. Bist, H. Oza, "Application of Machine Learning and Artificial Intelligence in Oil and Gas Industry", Petroleum Research, Vol. 6, No. 4, pp. 379-391, China, December 2021.
[2] A. Das, M.P. Sarma, K.K. Sarma, N. Mastorakis, "Design of an IoT Based Real Time Environment Monitoring System Using Legacy Sensors", Web of Conferences, Vol. 210, pp. 415-421, Majorca, Spain, July 2018.
[3] E.N. Aba, O.A. Olugboji, A. Nasir, M.A. Olutoye, O. Adedipe, "Petroleum Pipeline Monitoring Using an Internet of Things (IoT) Platform", SN Applied Sciences, Vol. 14, No. 2, pp. 1-12, January 2021.
[4] E.K. Ejomarie, E.G. Saturday, "Optimal Design of Gas Pipeline Transmission Network", GSJ, Vol. 8, No. 5, pp. 918-933, Rivers State, Nigeria, May 2020.
[5] F.I. Takhmazov, I.G. Aliyev, M.Z. Yusifov, N.A. Mammadova, "Technological Foundations of Multiline Gas Pipeline Reconstruction", International Journal on Technical and Physical Problems of Engineering (IJTPE), Issue 59, Vol. 16, No. 2, pp. 122-127, June 2024.
[6] I.G. Aliyev, "Working out the Method of Calculating the Dynamics Resulting from the Inclusion and Discharge of Gas Slug in the Operation of the Main Gas Pipeline", Eastern European Scientific Journal, Vol. 1, No. 100, pp. 25-34, Poland, May 2024.
[7] I.G. Aliyev, "Analysis of Transitional Processes in Various Operating Modes of Complex Gas Pipeline Systems", International Scientific and Practical Journal, No. 3, pp. 212-219, Kazakhstan, March 2024.
[8] J.N. Aslanov, K.S. Mammadov, "Design and Performance Analysis of Improved Valve Construction Being Used in Oil and Gas Industry", International Journal on Technical and Physical Problems of Engineering (IJTPE), Issue 51, Vol. 14, No. 2, pp. 98-103, June 2022.
[9] L.M.F. Silva, A.C.R. Oliveira, M.S.A. Leite, F.A.S. Marins, "Risk Assessment Model Using Conditional Probability and Simulation: Case Study in a Piped Gas Supply Chain in Brazil", International Journal of Production Research, Vol. 59, No. 10, pp. 2960-2976, England, May 2020.
[10] L. Meng, C. Liu, L. Fang, Y. Li, J. Fu, "Leak Detection of Gas Pipelines Based on Characteristics of Acoustic Leakage and Interfering Signals", SV Journals, Vol. 53, No. 4, China, 2019.
[11] M.A. Al Sabaeei, H. Alhussian, S.J. Abdulkadir, A. Jagadeesh, "Prediction of Oil and Gas Pipeline Failures Through Machine Learning Approaches: A Systematic Review", Energy Reports, No. 10, pp. 1313-1338, China, 2023.
[12] M.A. Adegboye, W.K. Fung, A. Karnik, "Recent Advances in Pipeline Monitoring and Oil Leakage Detection Technologies: Principles and Approaches", Sensors, Vol. 19, Article No. 11, pp. 25-48, UK, May 2019.
[13] M.R. Shadmesgaran, A.M. Hashimov, N.R. Rahmanov, "A Glance of Optimal Control Effects on Technical and Economic Operation in Grid", International Journal on Technical and Physical Problems of Engineering (IJTPE), Issue 46, Vol. 13, No. 1, pp. 1-10, March 2021.
[14] N.C. Ohalete, A.O. Aderibigbe, E.C. Ani, P.E. Ohenhen, A. Akinoso, "Advancements in Predictive Maintenance in the Oil and Gas Industry: A Review of AI and Data Science Applications", World Journal of Advanced Research and Reviews, Vol. 20, No. 3, pp. 167-181, Georgia, United States, 2023.
[15] N.E. Mohammad, Y.R. Yassmen, S. Aly, M.H. Hussein, "A Multi-Objective Optimization Method for Simulating the Operation of Natural Gas Transport System", Journal of Chemical Engineering, Vol. 41, No. 6, pp. 1609-1624, Korea, April 2024.
[16] O.A. Kurasov, P.V. Burkov, "Substantiation of Methods of Improving Safety of Pipeline Gas Transportation", E3S Web of Conferences, Vol. 266, pp. 01-012, Tomsk, 2021.
[17] R. Rios Mercado, Z. Borraz Sanchez, "Optimization Problems in Natural Gas Transportation Systems: A State-of-the-Art Review", Applied Energy, Elsevier, Vol. 147, pp. 536-555, Conrado, 2015.
[18] S. Li, C. Xia, Z. Cheng, W. Mao, Y. Liu, D. Habibi, "Leak Location Based on PDS-VMD of Leakage-Induced Vibration Signal Under Low SNR in Water-Supply Pipelines", IEEE Access, Vol. 8, pp. 68091-68102, China, March 2020.
[19] Y. Ding, P. Xu, Y. Lu, M. Yang, "Research on Pipeline Leakage Calculation and Correction Method Based on Numerical Calculation Method", Energies, Vol. 16, No. 21, p. 7255, China, October 2023.
[20] Z. Hafsi, A. Ekhtiari, L. Ayed, S. Elaoud, "The Linearization Method for Transient Gas Flows in Pipeline Systems Revisited: Capabilities and Limitations of the Modelling Approach", Journal of Natural Gas Science and Engineering, Vol. 101, pp. 104-494, Netherlands, 2022.
[21] Z. Li, Y. Liang, Y. Liang, Q. Liao, B. Wang, L. Huang, J. Zheng, H. Zhang, "Review on Intelligent Pipeline Technologies: A Life Cycle Perspective", Computers and Chemical Engineering, Vol. 175, pp. 108-283, China, July 2023.
[22] "Industrial Gas Pipeline Integrity Management", EIGA, Doc 235/21, p. 5, American, 2021.

## BIOGRAPHY

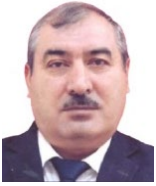

<u>Name</u>: **Ilgar**
<u>Middle Name</u>: **Giyas**
<u>Surname</u>: **Aliyev**
<u>Birthday</u>: 15.05.1966
<u>Birthplace</u>: Amasia, Armenia
<u>Education</u>: Construction Heating Gas Supply and Ventilation, Plumbing, Azerbaijan Architecture and Construction University, Baku, Azerbaijan, 1988
<u>Doctorate</u>: Technique/Construction and Operation of Oil and Gas Pipelines, Azerbaijan Architecture and Construction University, Baku, Azerbaijan, 1994
<u>The Last Scientific Position</u>: Assoc. Prof., Faculty of Civil Engineering, Azerbaijan University of Architecture and Construction, Baku, Azerbaijan, 1997
<u>Research Interests</u>: Construction and Operation of Oil-Gas Pipelines, Operation and Reconstruction of Buildings and Facilities, Designing Buildings and Facilities
<u>Scientific Publications</u>: 54, Papers 11, Books 16, 1 Patent, 12 Projects, 15 Theses